\documentclass[12pt]{amsart}

\usepackage{amsmath,amssymb,epic,lscape,lineno,setspace}

\newtheorem{theorem}{Theorem}[section]
\newtheorem{proposition}[theorem]{Proposition}

\theoremstyle{definition}

\theoremstyle{remark}

\numberwithin{equation}{section}
\providecommand{\bysame}{\leavevmode\hbox to3em{\hrulefill}\thinspace}

\def\DJ{{\hbox{D\kern-.8em\raise.15ex\hbox{--}\kern.35em}}}
\def\DJo{$\;$\kern-.4em
    \hbox{D\kern-.8em\raise.15ex\hbox{--}\kern.35em okovi\'c}}

\def\al{{\alpha}}

\def\la{{\lambda}}

\def\bZ{{\mbox{\bf Z}}}

\renewcommand{\subjclassname}{\textup{2000} Mathematics Subject
Classification }

\begin{document}

\title[Skew-Hadamard matrices]
{Skew-Hadamard matrices of orders $188$ and $388$ exist}

\author[D.\v{Z}. \DJ okovi\'{c}]
{Dragomir \v{Z}. \DJ okovi\'{c}}

\address{Department of Pure Mathematics, University of Waterloo,
Waterloo, Ontario, N2L 3G1, Canada}

\email{djokovic@uwaterloo.ca}

\thanks{
The author was supported by an NSERC Discovery Grant.}


\date{}

\begin{abstract}
We construct several difference families on cyclic groups of
orders $47$ and $97$, and use them to construct
skew-Hadamard matrices of orders $188$ and $388$.
Such difference families and matrices are constructed here
for the first time. The matrices are constructed by using
the Goethals--Seidel array.
\end{abstract}

\maketitle
\subjclassname{ 05B20, 05B30 }
\vskip5mm

\section{Introduction}

Recall that a Hadamard matrix $A$ of order $m$ is a $\{\pm1\}$-matrix
of size $m\times m$ such that $AA^T=mI_m$, where
$T$ denotes the transpose and $I_m$ the identity matrix.
A skew-Hadamard matrix is a Hadamard matrix $A$ such that
$A-I_m$ is a skew-symmetric matrix.
We refer the reader to \cite{CD} for the survey of
known results about skew-Hadamard matrices. 

The construction of skew-Hadamard matrices is lagging
considerably behind that for arbitrary Hadamard matrices.
Our previous four notes, written more than 13 years ago,
were motivated by the desire to improve this situation.
We constructed skew-Hadamard matrices of order
$m=4n$ for the following 24 odd integers $n$:
\begin{center}
\[ 
\begin{array}{l}
\text{\cite{DZ1}:} \quad 37,43; \\
\text{\cite{DZ2}:} \quad 67,113,127,157,163,181,241; \\
\text{\cite{DZ3}:} \quad 39,49,65,93,121,129,133,217,219,267; \\
\text{\cite{DZ5}:} \quad 81,103,151,169,463.
\end{array}
\]
\end{center}
At the time of publication, such matrices of these orders were
not known to exist. Due to the manifold increase in computing
power since that time, one can now make further progress.

In \cite{DZ5}, we listed 45 odd integers $n<300$ for which no
skew-Hadamard matrix of order $4n$ was known at that time.
(In the first edition of \cite{CD}, Table 24.31 was incomplete.)
The smallest of these $n$'s was $47$.
The next one, $59$, has been removed recently by Fletcher,
Koukouvinos and Seberry \cite{FKS}.
In this note we shall remove the integers $47$ and $97$ from the
mentioned list by constructing examples of skew-Hadamard
matrices of orders $4\cdot47=188$ and $4\cdot97=388$.
(We have constructed a bunch of examples but we have saved
and will present only a few of them.)

Consequently, the revised list now consists of the 42 integers:
\begin{eqnarray*}
&& 69,89,101,107,109,119,145,149,153,167,177,179,191,193, \\
&& 201,205,209,213,223,225,229,233,235,239,245,247,249,251, \\
&& 253,257,259,261,265,269,275,277,283,285,287,289,295,299.
\end{eqnarray*}

We construct our examples of skew-Hadamard matrices of
orders 188 and 388 by constructing first suitable supplementary
difference sets, and then we use these sets to build four
circulant blocks, which one should plug into the Goethals--Seidel
array. The procedure used to find these supplementary difference 
sets is not new. I have used it in several papers during the 
last 15 years. It is described in my note \cite{DZ4}.

\section{The case $n=47$}

We denote the additive group of integers modulo $n$ by $\bZ_n$.
In this section we set $n=47$. In the literature on Hadamard 
matrices it is customary to refer to difference families (DF)
as supplementary difference sets (SDS) and to employ
more elaborate and more informative notation by listing the
order $v$ of the underlying abelian group, the number of
sets in the family as well as their cardinals, and also the
parameter $\la$.

We have constructed four suitable difference
families in $\bZ_n$. The first two are the following.

\begin{proposition} \label{stav1}
Define six subsets of $\bZ_{47}$:
\begin{eqnarray*}
X_1 &=& \{ 2,3,5,6,7,9,10,11,12,13,14,17,18,19,20,21,22,25,27,30, \\
&& \quad	31,33,35,37,38,39,40,42,43,44 \}, \\
X_2 &=& \{ 1,3,6,7,8,11,13,14,15,19,20,21,24,27,30,33,39,41,43, \\
&& \quad 44,45,46 \}, \\
X_3 &=& \{ 3,6,8,10,11,12,14,20,21,23,24,25,26,27,30,31,32,34,35, \\
&& \quad 41,42,45 \}, \\
Y_1 &=& \{ 1,2,3,4,5,6,10,11,12,13,14,15,17,18,19,21,23,24,25,27,  \\
&& \quad	28,29,30,31,35,38,41,43,44,46 \}, \\
Y_2 &=& \{ 3,6,7,8,10,11,12,16,22,25,26,31,32,33,34,37,39,41,42, \\
&& \quad 43,44,46 \}, \\
Y_3 &=& \{ 3,7,12,13,15,16,18,20,21,23,25,26,27,28,32,35,38,39,42, \\
&& \quad 44,45,46 \}.
\end{eqnarray*}
The triples $\{X_1,X_2,X_3\}$ and $\{Y_1,Y_2,Y_3\}$ are
difference families, i.e., they are $3-(47;30,22,22;39)$
supplementary difference sets in $\bZ_{47}$.
The two families are not equivalent.
\end{proposition}
\begin{proof} Use the computer to verify the claims. Note that
the cardinals $n_k=|X_k|=|Y_k|$ are indeed $n_1=30$ and $n_2=n_3=22$.
The parameter $\la$ is $39$, i.e., each nonzero integer in
$\bZ_n$ occurs $39$ times in the list of differences created
from the sets $X_k$ and also from the $Y_k$.

The second claim can be verified in several ways. We used the
following ad hoc method.
We compare the list of differences generated by
the sets $X_1$ and $Y_1$. Each nonzero integer $i\in\bZ_n$ occurs
in one of these lists say $\mu_i$ times. The $\mu_i$'s take only
three values: 18, 19 or 20. But the number of $\mu_i$'s equal to
18, 19 and 20 is 12, 26 and 8 for $X_1$ and 14, 22 and 10 for $Y_1$.
Hence $X_1$ and $Y_1$ are not equivalent under translations and
automorphisms of the additive group $\bZ_n$.
\end{proof}

For any subset $X\subseteq\bZ_n$ let
\[ a_X=(a_0,a_1,\ldots,a_{n-1}) \]
be the $\{\pm1\}$-row vector such that $a_i=-1$ iff $i\in X$.
We denote by $A_X$ the $n\times n$ circulant matrix having
$a_X$ as its first row.

Let $X_0\subseteq\bZ_n$ be the Paley difference set (the set of
nonzero squares in the finite field $\bZ_n$). Recall that $X_0$
is of skew type, i.e., for nonzero $i\in\bZ_n$ we have
$i\in X_0$ iff $-i\notin X_0$. Its cardinal is $n_0=|X_0|=23$.

For simplicity, write $A_k$ instead of $A_{X_k}$ for $k=0,1,2,3$.
We can now plug our matrices $A_k$ into the Goethals--Seidel template
to construct a skew-Hadamard matrix of order $188$:
\[
A=\left[ \begin{array}{cccc}
		A_0 & A_1R & A_2R & A_3R \\
		-A_1R & A_0 & -A_3^TR & A_2^TR \\
		-A_2R & A_3^TR & A_0 & -A_1^TR \\
		-A_3R & -A_2^TR & A_1^TR & A_0
\end{array} \right]. \]
As usual, $R$ denotes the matrix having ones on the back-diagonal
and all other entries zero.
		
Clearly, we can use the second difference family to construct
another skew-Hadamard matrix of order 188. Both solutions have
the same associated decomposition of $4n$ as sum of four squares:
\begin{eqnarray*}
4n=188 &=& 13^2+3^2+3^2+1^2 \\
&=& \sum_{k=0}^3 (n-2n_k)^2.
\end{eqnarray*}

The remaining two difference families have different parameters
from the first two.

\begin{proposition} \label{stav2}
Define six subsets of $\bZ_{47}$:
\begin{eqnarray*}
P_1 &=& \{ 0,2,4,5,9,10,12,16,17,19,21,22,23,25,27,28,35,36,37,  \\
&& \quad	43,46 \}, \\
P_2 &=& \{ 0,1,2,6,8,9,11,15,16,19,25,32,33,35,36,37,38,40,44 \}, \\
P_3 &=& \{ 1,2,3,4,5,6,7,10,11,16,18,22,24,28,31,35,38,40,43 \}, \\
Q_1 &=& \{ 4,5,6,8,11,12,15,20,21,23,25,26,28,29,30,31,32,36, \\
&& \quad	39,41,43 \}, \\
Q_2 &=& \{ 1,2,5,7,13,14,21,22,24,26,31,32,35,36,37,39,40,42,46 \}, \\
Q_3 &=& \{ 1,2,3,4,5,9,12,18,20,21,24,25,32,34,38,39,43,44,46 \}.
\end{eqnarray*}
The triples $\{P_1,P_2,P_3\}$ and $\{Q_1,Q_2,Q_3\}$ are
difference families, i.e., they are $3-(47;21,19,19;24)$
supplementary difference sets in $\bZ_{47}$. These two families are not
equivalent to each other or the ones above.
\end{proposition}

Just as the first two families, $\{P_1,P_2,P_3\}$ and 
$\{Q_1,Q_2,Q_3\}$ can be used to construct two more skew-Hadamard 
matrices of order 188. The associated decomposition into
sum of four squares is now different: $188 = 9^2+9^2+5^2+1^2.$

\section{The case $n=97$}

For the remainder of this note we set $n=97$.
Let $G$ be the multiplicative group of the nonzero elements of
$\bZ_n$, a cyclic group of order $n-1=96$, and let
$H=\langle 35 \rangle =\{ 1,35,61 \}$ be its subgroup of order 3.
We use the same enumeration of the 32 cosets $\al_i$,
$0\le i\le 31$, of $H$ in $G$ as in our computer program.
Thus we impose the condition that $\al_{2i+1}=-1\cdot\al_{2i}$
for $0\le i\le 15$. For even indices we have
\[
\begin{array}{lllll}
\al_0=H, \quad & \al_2=2H, \quad & \al_4=3H, \quad & \al_6=4H, \quad & \al_8=5H, \\
\al_{10}=6H, & \al_{12}=7H, & \al_{14}=9H, & \al_{16}=10H, & \al_{18}=12H, \\
\al_{20}=13H, & \al_{22}=15H, & \al_{24}=18H, & \al_{26}=20H, & \al_{28}=23H, \\
\al_{30}=26H. &&&&
\end{array}
\]

Next define four index sets:
\begin{eqnarray*}
J_0 &=& \{ 1,2,4,6,9,11,13,14,17,18,21,23,25,27,29,30 \}, \\
J_1 &=& \{ 1,2,6,7,8,9,10,11,12,13,23,27,29 \}, \\
J_2 &=& \{ 0,1,2,5,6,12,13,15,16,20,24,25,26,29,30,31 \}, \\
J_3 &=& \{ 0,2,3,4,7,8,9,11,12,13,15,16,17,18,23,28,29 \}
\end{eqnarray*}
and introduce the following four subsets of $\bZ_n$:
\[ U_k = \bigcup_{i\in J_k} \al_i,\quad k=0,1,2,3. \]
Their cardinals $n_k=|U_k|=3|J_k|$ are:
\[ n_0=n_2=48,\, n_1=39,\, n_3=51 \]
and we set
\[ \la = n_0+n_1+n_2+n_3-n=89. \]
Observe that $U_0$ is of skew type, i.e., we have
\[ U_0\cap (-U_0)=\emptyset,\quad  U_0\cup (-U_0)=\bZ_n\setminus\{0\}. \]

\begin{proposition} \label{stav3}
The four subsets $U_0,U_1,U_2,U_3\subset\bZ_n$
form a difference family, i.e., they are $4-(97;48,39,48,51;89)$
supplementary difference sets in $\bZ_{97}$.
\end{proposition}
\begin{proof}
For $r\in\{ 1,2,\ldots,96 \}$ let $\la_k(r)$ denote the number of
solutions of the congruence $i-j\equiv r \pmod{97}$ with
$\{i,j\}\subseteq U_k$. It is easy to verify (by using a computer) that
\[ \la_1(r)+\la_2(r)+\la_3(r)+\la_4(r)=\la \]
is valid for all such $r$. Hence the sets $U_1,U_2,U_3,U_4$
form a difference family in $\bZ_n$.
\end{proof}

Let $A_k$ now denote the $n\times n$ circulant matrices $A_{Y_k}$.
The SDS-property implies that the $\{\pm1\}$-matrices
$A_0,\ldots,A_3$ satisfy the identity
\[ \sum_{k=0}^3 A_kA_k^T=4nI_n. \]

One can now plug the matrices $A_k$ into the Goethals--Seidel
template to obtain a Hadamard matrix $A$ of order $4n=388$.
Since $U_1$ is of skew type, $A$ is also skew-Hadamard.

Our second example, $B$, is constructed in the same way by using
the index sets:
\begin{eqnarray*}
K_0 &=& \{ 0,3,4,7,9,11,12,14,17,19,20,22,24,27,28,30 \}, \\
K_1 &=& \{ 4,7,8,10,12,13,14,15,17,18,20,26,27 \}, \\
K_2 &=& \{ 0,1,2,3,6,7,8,11,12,14,20,23,24,25,28,31 \}, \\
K_3 &=& \{ 1,2,4,7,8,9,10,12,13,19,21,23,24,25,26,27,31 \},
\end{eqnarray*}
with the corresponding subsets of $\bZ_n$:
\[ V_k = \bigcup_{i\in K_k} \al_i,\quad k=0,1,2,3, \]
with $V_0$ of skew type.

\begin{proposition} \label{stav4}
The four subsets $V_0,V_1,V_2,V_3\subset\bZ_n$
form a difference family, i.e., they are $4-(97;48,39,48,51;89)$
supplementary difference sets in $\bZ_{97}$.
\end{proposition}

The two SDS's that we used to construct $A$ and $B$ are not equivalent.
For instance, the sets $U_1$ and $V_1$ are not equivalent
under translations and group automorphisms of $\bZ_n$.

Since the two SDS's have the same parameters, they share
the same decomposition of $4n$ into sum of four squares:
\begin{eqnarray*}
4n=388 &=& 19^2+5^2+1^2+1^2 \\
&=& \sum_{k=1}^4 (n-2n_k)^2.
\end{eqnarray*}

\end{document}